\documentclass[11pt,a4paper]{amsart}
\usepackage{amsmath,amsfonts,amssymb,epsf,mathrsfs,xypic}
\ifx\optionkeymacros\undefined\else \fi

\catcode`\å=\active\defå{{\aa}}       % option a
\catcode`\†=\active\def†{\int}        % option b (math mode) 
\catcode`\ç=\active\defç{\c c}        % option c
\catcode`\=\active\def{\partial}    % option d (math mode)
\catcode`\Ÿ=\active\defŸ{\oint}       % option f (math mode) ?
\catcode`\=\active\def{\triangle}   % option j (math mode)
\catcode`\¬=\active\def¬{\neg}        % option l (math mode)
\catcode`\µ=\active\defµ{\mu}         % option m (math mode)
\catcode`\ø=\active\defø{{\o}}        % option o
\catcode`\¼=\active\def¼{\pi}         % option p (math mode w/ arg.)
\catcode`\¦=\active\def¦{{\oe}}       % option q 
\catcode`\ß=\active\defß{{\ss}}       % option s 
\catcode`\Ý=\active\defÝ{\dagger}     % option t  (math mode)
\catcode`\ˆ=\active\defˆ{\sqrt}       % option v (math mode w/ arg.)
\catcode`\…=\active\def…{\Sigma}      % option w (math mode)
\catcode`\‰=\active\def‰{\approx}     % option x (math mode)
\catcode`\‡=\active\def‡{\Omega}      % option z (math mode)
\catcode`\£=\active\def£{{\it\$}}     % option 3 ($ from italic font)
\catcode`\ƒ=\active\defƒ{\infty}      % option 5 (math mode)
\catcode`\§=\active\def§{{\S}}        % option 6 
\catcode`\¶=\active\def¶{{\P}}        % option 7
\catcode`\€=\active\def€{\bullet}     % option 8 
\catcode`\ª=\active\defª{\leavevmode\raise.585ex\hbox{\b a}}      % option 9
\catcode`\º=\active\defº{\leavevmode\raise.6ex\hbox{\b o}}        % option 0
\catcode`\'=\active\def'{\not=}       % option = (math mode)
\catcode`\¾=\active\def¾{\leq}        % option , (math mode)
\catcode`\"=\active\def"{\geq}        % option . (math mode)
\catcode`\÷=\active\def÷{\div}        % option / (math mode)
\catcode`\Š=\active\defŠ{{\dots}}     % option ; 
\catcode`\æ=\active\defæ{{\ae}}       % option '
\catcode`\«=\active\def«{\ll}         % option \ (math mode)
\catcode`\³=\active\def³{``}          % option [
\catcode`\¡=\active\def¡{!`}          % option !
\catcode`\¢=\active\def¢{\rlap/c}     % option 4
\catcode`\Œ=\active\defŒ{`}           % option ] 
\catcode`\¹=\active\def¹{'}           % shift option ]

\catcode`\Å=\active\defÅ{{\AA}}       % shift-option A
\catcode`\Ç=\active\defÇ{\c C}        % shift-option C
\catcode`\Ø=\active\defØ{{\O}}        % shift-option O
\catcode`\½=\active\def½{\Pi}         % shift-option P (math mode)
\catcode`\'=\active\def'{{\OE}}       % shift-option Q
\catcode`\Æ=\active\defÆ{{\AE}}       % shift-option '
\catcode`\×=\active\def×{\diamond}    % shift-option V (math mode)
\catcode`\°=\active\def°{\accent'27}  % shift-option 8
\catcode`\²=\active\def²{''}          % shift-option [
\catcode`\±=\active\def±{\pm}         % shift-option = (math mode)
\catcode`\»=\active\def»{\gg}         % shift-option \ (math mode)
\catcode`\¿=\active\def¿{?`}          % shift-option / 
\catcode`\­=\active\def­{--}          % option - (en-dash)
\catcode`\‹=\active\def‹{---}         % shift-option - (em-dash)

\catcode`\ä=\active\defä{\"a}        % option u, then  a
\catcode`\ë=\active\defë{\"e}        % option u, then  e
\catcode`\ï=\active\defï{\"{\i}}     % option u, then  i
\catcode`\ö=\active\defö{\"o}        % option u, then  o
\catcode`\ü=\active\defü{\"u}        % option u, then  u
\catcode`\ÿ=\active\defÿ{\"y}        % option u, then  y
\catcode`\Ä=\active\defÄ{\"A}        % option u, then  A
\catcode`\Ö=\active\defÖ{\"O}        % option u, then  O
\catcode`\Ü=\active\defÜ{\"U}        % option u, then  U
\catcode`\á=\active\defá{\'a}        % option e, then  a
\catcode`\é=\active\defé{\'e}        % option e, then  e
\catcode`\í=\active\defí{\'{\i}}     % option e, then  i
\catcode`\ó=\active\defó{\'o}        % option e, then  o
\catcode`\ú=\active\defú{\'u}        % option e, then  u
\catcode`\É=\active\defÉ{\'E}        % option e, then  E
\catcode`\à=\active\defà{\`a}        % option `, then  a
\catcode`\è=\active\defè{\`e}        % option `, then  e
\catcode`\ì=\active\defì{\`{\i}}     % option `, then  i
\catcode`\ò=\active\defò{\`o}        % option `, then  o
\catcode`\ù=\active\defù{\`u}        % option `, then  u
\catcode`\À=\active\defÀ{\`A}        % option `, then  A
\catcode`\ã=\active\defã{\~a}        % option n, then  a
\catcode`\ñ=\active\defñ{\~n}        % option n, then  n
\catcode`\õ=\active\defõ{\~o}        % option n, then  o
\catcode`\Ã=\active\defÃ{\~A}        % option n, then  A
\catcode`\Ñ=\active\defÑ{\~N}        % option n, then  N
\catcode`\Õ=\active\defÕ{\~O}        % option n, then  O
\catcode`\â=\active\defâ{\^a}        % option i, then  a
\catcode`\ê=\active\defê{\^e}        % option i, then  e
\catcode`\î=\active\defî{\^{\i}}     % option i, then  i
\catcode`\ô=\active\defô{\^o}        % option i, then  o
\catcode`\û=\active\defû{\^u}        % option i, then  u

\let\optionkeymacros\null

\usepackage[backref=page,bookmarksopen=true]{hyperref}

\let\bbbibitem\bibitem
\renewcommand{\bibitem}[2][]{\bbbibitem[#1]{#2}\label{#2}}

\def\fin{\hfill\hbox{\hskip .2cm $\square$}\medskip}%{\vrule height 5pt width 5 pt}}}

\newtheorem{theo}{Theorem}[section]
\newtheorem{lemma}[theo]{Lemma}
\newtheorem{prop}[theo]{Proposition}

\newtheorem{defi}[theo]{Definition}
\newtheorem{propdef}[theo]{Proposition-definition}

\newenvironment{demo}[1][\hspace{-3pt}]{{\noindent\em Proof #1.~ }}{\fin}

\newenvironment{thm}[1][\hspace{-3pt}]{{\vspace{0.4cm}\noindent\bf 
Theorem #1}~~\em}{\vspace{0.4cm}}

\def\cal{\mathcal}
\def\a{\alpha}
\def\b{\beta}
\def\deg{{\rm deg}\,}
\def\g{\gamma}
\def\G{\Gamma}
\def\R{{\mathbb R}} 
\def\d{{\rm d}}

\def\ad{{\rm ad}}

\def\PU{{\rm PU}}
\def\SU{{\rm SU}}
\def\GL{{\rm GL}}

\def\Hom{{\rm Hom}}
\def\U{{\rm U}}
\def\SL{{\rm SL}}
\def\SO{{\rm SO}}

\def\Im{{\rm Im\,}}
\def\C{{\mathbb C}}
\def\Q{{\mathbb Q}}
\def\B{{\mathbb B}}
\def\Bm{{\mathbb B}^m}    

\def\N{{\mathbb N}}
\def\Z{{\mathbb Z}}
\def\R{{\mathbb R}}
\def\fd{\longrightarrow}

\def\rk{{\rm rk}\,}
\def\ca{{\rm can}}

\def\df{{\rm d}f}

\def\ds{\displaystyle}

\def\wt{\widetilde}

\def\o{\omega}

\def\X{{\cal X}}
\def\Y{{\cal Y}}
\def\U{{\rm U}}

\def\L{{\cal L}}

\def\gg{{\mathfrak g}}
\def\kg{{\mathfrak k}}

\def\pg{{\mathfrak p}}
\def\hg{{\mathfrak h}}
\def\qg{{\mathfrak q}}
\def\ug{{\mathfrak u}}

\def\zg{{\mathfrak z}}

\def\P{{\mathbb P}}

\def\tr{{\rm tr}\,}
\def\trans{{}^t\hspace{-2pt}}

\def\V{{\mathbb V}}
\def\W{{\mathbb W}}
\def\E{{\mathbb E}}
\begin{document}

\title[The Toledo invariant on smooth varieties of general type]{The Toledo invariant\\on smooth varieties of general type}
%\titlerunning{Representations of uniform complex hyperbolic lattices}

\author{Vincent Koziarz} 
\author{Julien Maubon}
\address{IECN, Nancy-Université, CNRS, INRIA, Boulevard des Aiguillettes B. P. 239, F-54506 
Vand\oe uvre-lès-Nancy, France}
\email{koziarz@iecn.u-nancy.fr}  
\email{maubon@iecn.u-nancy.fr}

%\date{\today}

\sloppy

\begin{abstract}  
We propose a definition of the Toledo invariant for representations of fundamental groups of smooth varieties of general type into semisimple Lie groups of Hermitian type. This definition allows to generalize the results known in the classical case of representations of complex hyperbolic lattices to this new setting: assuming that the rank of the target Lie group is not greater than two, we prove that the Toledo invariant satisfies a Milnor-Wood type inequality and we characterize the corresponding maximal representations.      
\end{abstract}

\maketitle

\section{Introduction}

The Toledo invariant is a characteristic number naturally associated to representations of lattices of semisimple Lie groups of Hermitian type into other semisimple Lie groups of Hermitian type. Recall that a real semisimple Lie group $G$ (with no compact factors) is said to be of Hermitian type if its associated symmetric space is Hermitian symmetric.   
The most general definition of this invariant can be found in Burger and Iozzi paper~\cite{BI}. It is given there in terms of the second bounded cohomology of the involved Lie groups and lattices but, in the case of torsion free uniform lattices, we can rephrase it as follows. Let $\X$ be an irreducible Hermitian symmetric space of the noncompact type and let $\Gamma$ be a (torsion free) uniform lattice in the automorphism group of $\X$. In this introduction, we loosely speak of the ``automorphism group of $\X$'' while we actually mean ``a Lie group acting transitively and almost effectively by complex automorphisms on $\X$'', which could therefore be a finite cover of the true automorphism group of $\X$. Let $\rho$ be a representation, i.e. a group homomorphism, of $\G$ into a linear connected simple noncompact Lie group of Hermitian type $G$. Let also $\Y$ be the Hermitian symmetric space associated to $G$, so that $\Y=G/K$ for some maximal compact subgroup $K$ of $G$.
There always exists a smooth $\rho$-equivariant map $f:\X\fd\Y$ and we can use this map to pull-back the $G$-invariant Kähler form $\o_\Y$ of $\Y$ to a closed 2-form $f^\star\o_\Y$ on $\X$ which by equivariance goes down to a form on $X:=\G\backslash\X$ that will still be denoted by $f^\star\o_\Y$. The cohomology class defined by $f^\star\o_\Y$ depends only on $\rho$, not on the chosen equivariant map $f$. This class can then be evaluated against the Kähler class $\o_X$ of $X$ coming from the invariant Kähler form $\o_\X$ on $\X$ and this gives the Toledo invariant of $\rho$:     
$$
\tau(\rho):=\frac{1}{m!}\int_X f^\star\o_\Y\wedge\o_X^{m-1},
$$
where $m$ is the complex dimension of $X$.

Burger and Iozzi proved that the Toledo invariant satisfies a Milnor-Wood type inequality, namely that it is bounded in absolute value by a constant depending only on the ranks of the symmetric spaces $\X$ and $\Y$ and the volume of $X=\G\backslash \X$: if the invariant Kähler metrics $\o_\X$ and $\o_\Y$ are normalized so that the minimum of their holomorphic sectional curvatures is $-1$, then 
$$
|\tau(\rho)|\leq\frac{\rk\Y}{\rk\X}{\rm Vol}(X)
$$
where ${\rm Vol}(X)$ is computed w.r.t. $\o_X$, i.e. ${\rm Vol}(X)=\frac{1}{m!}\int_X \o_X^{m}$.

One of the main motivations in studying the Toledo invariant is that it allows to prove rigidity results by singling out a special class of representations, those for which the Milnor-Wood inequality is an equality. These so-called maximal representations have received a great amount of attention and one expects to be able to describe them completely. In this rigidity oriented approach, the interesting case is the case of lattices $\G$ in rank one Lie group (by the rank of a semisimple Lie group we shall always mean its split rank, that is the rank of its associated symmetric space, so that for example $\rk G=\rk\Y$). Indeed, lattices in higher rank Lie groups are known to be superrigid, which is stronger than what can be proved using the Toledo invariant. Therefore in this paper $\X$ will always be complex hyperbolic $m$-space, which as a bounded symmetric domain is the unit ball $\Bm$ in $\C^m$, and $\G$ will be a (torsion free) uniform lattice in the automorphism group of $\Bm$. Such a lattice will be called a surface group if $m=1$, that is when $X=\G\backslash\B^1$ is a Riemann surface, and a higher dimensional complex hyperbolic lattice if $m\geq 2$.       

The Toledo invariant was first defined for representations of surface groups into $G=\SU(n,1)$ in 1979 Toledo's paper~\cite{T1} and then more explicitly in~\cite{T2}, where the Milnor-Wood inequality was proved in this case. Toledo also proved that maximal representations are faithful with discrete image, and stabilize a complex line in the complex hyperbolic $n$-space $\Y$. At approximately the same time, Corlette established in~\cite{C} the same kind of result for higher dimensional lattices (he used a different but very similar invariant), and this was extended to non uniform lattices in~\cite{BI} and~\cite{KM0} using the Toledo invariant. An immediate corollary is that a lattice in $\SU(m,1)$ can not be deformed non-trivially in $\SU(n,1)$, $n\geq m\geq 2$, a result first obtained in the uniform case in~\cite{GM} using different methods. Therefore the case where the rank of the Lie group $G$ is one is now settled.  Maximal representations of surface groups into higher rank Lie groups are also quite well understood, thanks to the work of Burger, Iozzi, Wienhard~\cite{BIW} and Bradlow, Garcia-Prada, Gothen~\cite{BGPG1,BGPG2}. Concerning higher dimensional lattices, the case of classical target Lie groups of rank two has been (almost entirely) treated in our previous paper~\cite{KM}, but no general proof of the expected rigidity of maximal representations into higher rank Lie groups is known.

\medskip

From a somewhat different perspective, the lattices we are considering are examples of Kähler groups, i.e. fundamental groups of closed Kähler manifolds. Being a Kähler group is a severe restriction and these groups share many rigidity properties with lattices in semisimple Lie groups, the first historical example of this phenomenon being Siu's strengthening of Mostow strong rigidity theorem~\cite{Si80}. The aim of this paper is to show how the definition of the Toledo invariant and the known characterizations of maximal representations can be generalized to the case where the represented group $\G$ is no more a (higher dimensional) complex hyperbolic lattice, but merely the fundamental group of a smooth variety of general type.  

This generalization is made possible by the fact that some of the results we mentioned were (or can be) proved with complex differential geometric methods. In particular, the most recent works use the theory of Higgs bundles, developed by Hitchin\cite{H1,H2} and Simpson~\cite{S1,S2} precisely for studying linear representations of Kähler groups (see \cite{BGPG1,BGPG2} for surface groups, and \cite{KM} for higher dimensional lattices). This approach casts a new light on the reason why such results hold and indicates that they should be valid in a much wider setting. 

Before explaining a little the issues one has to address when trying to generalize the definition of the Toledo invariant, let us say a few words about why we think the context of varieties of general type is the right one. First, there are many examples of interesting representations of fundamental groups of varieties of general type into noncompact Lie groups of Hermitian type. For instance, this context encompasses highly non trivial representations $\pi_1(X)\fd\SU(m,1)$ for some Kähler manifolds $X$ which have negative sectional curvature, but are not covered by the complex unit ball $\Bm$ of dimension $m$. We think of the so-called Mostow-Siu surfaces and their generalizations (see the more recent work of Deraux~\cite{Der1}, and~\cite{Der2} for a 3-dimensional example) which by construction admit such representations. On the other hand, although we are not dealing with general Kähler groups, 
a factorization theorem of Zuo~\cite{Z}, generalizing~\cite{Mo}, asserts that given a compact Kähler manifold $X$ and a Zariski-dense representation $\rho:\pi_1(X)\fd G\subset {\rm SL}(n,\R)$ of its fundamental group into a linear simple Lie group, there exists a finite etale cover $e:X^e\fd X$, a proper modification $\hat X^e\fd X^e$, a smooth variety $X_{gt}$ of general type, a surjective holomorphic map $s:\hat X^e\fd X_{gt}$ and a representation $\rho_{gt}:\pi_1(X_{gt})\fd G$ such that $\rho\circ e_\star=\rho_{gt}\circ s_\star$, where $\pi_1(\hat X^e)$ is canonically identified with $\pi_1(X^e)$. So, in a way, the study of  representations of fundamental groups of Kähler manifolds in linear simple Lie groups is reduced to the study of representations of fundamental groups of smooth varieties of general type.

For the complex differential geometric viewpoint we want to adopt, it is more convenient to consider the Toledo invariant as the degree of a line bundle over $X$ associated to the representation $\rho:\G=\pi_1(X)\fd G$. This interpretation is indeed the best suited for the generalization to representations of fundamental groups of varieties of general type we are aiming at. Even in the classical case, namely when studying representations of complex hyperbolic lattices, it is difficult to give definitions and proofs as independent as possible of the different types of target Lie groups $G$. In the end one almost always needs to resort to the classification of Hermitian symmetric spaces. We have included in Section~\ref{classical} a somehow detailed discussion of this question and the reader more interested by varieties of general type may read only the beginning, skip to Definition~\ref{defclassical}, and then to Section~\ref{general}, where we define the Toledo invariant for fundamental groups of varieties of general type. 
Whereas in the classical case degrees are naturally computed using the canonical polarization of the locally symmetric manifold $\G\backslash\Bm$, there is no such obvious choice of a polarization on a smooth variety of general type $X$. Examination of the existing proofs of the Milnor-Wood inequality shows in particular that in order to get a useful definition of the Toledo invariant in the general type case, namely one for which such an inequality holds, we need a polarization with respect to which the tangent bundle of $X$ is semistable. It turns out that in general it is necessary to choose what we call a good birational model of the variety whose fundamental group we are representing, for polarizations with the right properties do not exist on all birational models. 
The existence of these good models relies on the existence of the canonical model of the variety $X$, a recent and very deep result (\cite{BCHM}). 
Since any two smooth birational varieties have isomorphic fundamental groups, $\rho$ can be considered as a representation of both those groups. As a consequence, the represented group $\G$ should not be understood as the fundamental group of a precise smooth variety of general type $X$ but rather as the common fundamental group of all smooth representatives of the birational class of $X$, and the representation $\rho$ should be seen as a representation of this kind-of-abstract group. The definition of the Toledo invariant we propose indeed does not depend on $X$, only on its birational class.

With this in mind, we can extend the results of~\cite{KM} to our new setting. We consider representations into Lie groups whose associated symmetric spaces are the classical Hermitian symmetric spaces of noncompact type of rank one or two (with the exception of the symmetric space associated to $\SO^\star(10)$). This condition on the rank of the target Lie group comes from the fact that, already in the classical case of representations of higher dimensional complex hyperbolic lattices, the Milnor-Wood type inequality of Burger-Iozzi has not yet been proved within the complex geometric framework under more general assumptions. Such restrictions, on the ranks of the groups or equivalently on the ranks of the associated Higgs bundles, appear for similar reasons when dealing with numerous related questions, see for example~\cite{Kl} and~\cite{Rez}. On the other hand, the known complex differential geometric proofs of the Milnor-Wood inequality give more information on the equality case than Burger and Iozzi general proof, and indeed they allow to completely characterize maximal representations. Here, because of the nature of the Toledo invariant, the characterization we obtain concerns both the representation and the variety whose fundamental groups is maximally represented, or more precisely its birational class, through its canonical model.    

\begin{thm}[\hspace{-4pt}.]
Let $X$ be a smooth variety of general type and of complex dimension $m\geq 2$ and let $X_\ca$ be its canonical model. 
Let $G$ be either $\SU(p,q)$ with $1\leq q\leq 2\leq p$, ${\rm Spin}(p,2)$ with $p\geq 3$, or ${\rm Sp}(2,\R)$. 
Finally let $\rho:\pi_1(X)\fd G$ be a representation. 

Then the  Toledo invariant of $\rho$ satisfies the Milnor-Wood type inequality  
$$
|\tau(\rho)|\leq \rk G\,\frac{K_{X_\ca}^m}{m+1},
$$
where $K_{X_\ca}$ is the canonical divisor of $X_\ca$.  

Equality holds if and only if $G=\SU(p,q)$ with $p\geq qm$ and there exists a $\rho$-equivariant (anti)holomorphic proper embedding from the universal cover of $X_\ca$ onto a totally geodesic copy of complex hyperbolic $m$-space $\Bm$, of induced holomorphic sectional curvature $-1/q$, in the symmetric space associated to $G$. In particular, the canonical model of $X$ is then smooth and uniformized by $\Bm$, and the representation $\rho$ is discrete and faithful.
\end{thm}

The proof of the theorem is given in Section~\ref{proof}. In Subsection~\ref{inequality} we prove the Milnor-Wood type inequality. We begin by reviewing the stability properties, for vector bundles as well as for Higgs bundles, that were already needed in the classical case of representations of complex hyperbolic lattices, and we check that they still hold for the polarization chosen to define the Toledo invariant in the general type case, the crucial point being the semistability of the tangent bundle of the variety. Once these properties are established, the proof of the Milnor-Wood inequality goes exactly like in the classical case. Maximal representations are discussed in Subsection~\ref{equality}. Their characterization is noticeably more difficult than in the classical case, because the Ahlfors-Schwarz-Pick lemma is not available anymore, and because we need to work on the possibly singular canonical model of the variety.   

\medskip

{\em Acknowledgments.} We would like to thank Frédéric Campana for sharing with us his knowledge and experience in the course of  numerous discussions about the Minimal Model Program. We are grateful to Stéphane Druel for his explanations about the paper~\cite{BCHM}. We also benefited from useful conversations with Jean-Louis Clerc and with Ngaiming Mok, who we particularly thank for drawing our attention to his joint paper~\cite{EM} with Philippe Eyssidieux.  

\section{The Toledo invariant as the degree of a line bundle}\label{classical}

Given a nef $\Q$-line bundle $L$ (a polarization) on a variety $X$, the $L$-degree of any coherent sheaf ${\cal F}$ on $X$ is given by: 
$$
{\rm deg}_L\, {\cal F}:=c_1({\cal F})\cdot c_1(L)^{m-1}=\int_X c_1({\cal F})\wedge c_1(L)^{m-1}
$$
where $c_1({\cal F})\in H^2(X,\R)$ is the first Chern class of ${\cal F}$ (see for example~\cite{K} for the definition) and $c_1(L)\in H^2(X,\R)$ denotes the first Chern class of $L$. In the same way, we can compute the $L$-degree of any smooth complex line bundle $F$ and we will use freely the notations ${\rm deg}_L \,F=F\cdot L^{m-1}$. If moreover $F$ is associated to a divisor $D$, we will also write $D\cdot L^{m-1}$ for ${\rm deg}_L \,F$.

\medskip

When $X$ is a closed Hermitian locally symmetric space of the noncompact type and $\rho$ a representation of its fundamental group into a Lie group of Hermitian type, the Toledo invariant of $\rho$ can be interpreted as the degree of some line bundle on $X$ computed with the polarization coming from the canonical bundle $K_X$ of $X$. There is a certain degree of freedom in the choice of this line bundle and this is what we want to focus on in this section. Later, when dealing with fundamental groups of varieties of general type, it is the polarization used to compute degrees that we shall discuss.       

\medskip

What follows is valid with the obvious modifications for lattices in the automorphism group of any Hermitian symmetric space of the noncompact type but for the reasons given in the introduction we will stick here to complex hyperbolic lattices. So let $\G$ be a (torsion free) uniform lattice in the automorphism group of complex hyperbolic $m$-space $\Bm$, and let $X=\Gamma\backslash\Bm$ be the quotient complex hyperbolic manifold. Let $G$ be a linear connected simple noncompact Lie group of Hermitian type. Remember that $\rk G$ is the split rank of $G$, namely the rank of its associated symmetric space $\Y$, and that $\o_\Y$ and $\o_X$ are respectively the $G$-invariant Kähler form on $\Y$ and the Kähler form on $X$ coming from the invariant Kähler form on $\Bm$, both normalized so that the minimum of their holomorphic sectional curvatures is $-1$. Let finally $\rho$ be a representation of $\G$ in $G$.   

Hermitian symmetric spaces of the noncompact type are Kähler-Einstein manifolds with negative first Chern form, hence the first Chern form of their canonical bundle is their Kähler form up to a positive multiplicative constant. Therefore, if $c_\Y$ is such that $c_1(K_\Y)=\frac{c_\Y}{4\pi}\,\o_\Y$ and since $c_{\B^m}=m+1$, 
$$
\tau(\rho)=\frac{1}{m!}\frac{(4\pi)^m}{c_\Y\,(m+1)^{m-1}}\int_X c_1(f^\star K_\Y)\wedge c_1(K_X)^{m-1}=:\frac{1}{m!}\frac{(4\pi)^m}{c_\Y\,(m+1)^{m-1}}\deg(f^\star K_\Y),
$$
where degrees are computed w.r.t. the canonical polarization of $X$, $f:\Bm\fd \Y$ is any smooth $\rho$-equivariant map and $f^\star K_\Y$ is the vector bundle on $X$ obtained by first pulling-back the canonical bundle $K_\Y$ of $\Y$ to $\Bm$ by $f$ and then pushing it down to $X$. Note that the isomorphism class of $f^\star K_\Y$ depends only on $\rho$, not on $f$, since any two such equivariant maps are equivariantly homotopic. 
 
The Milnor-Wood inequality can then be written
$$
\left|\frac{\deg f^\star K_\Y}{c_\Y\,\rk \Y}\right|\leq \frac{\deg K_X}{m+1}.
$$

\medskip

We would like to get rid of the normalizing constant $c_\Y$ in order to make our forthcoming statements as independent as possible of the different types of target Lie groups or symmetric spaces. For this, we embed $\Y$ into its compact dual $\Y_c$, and we use the fact that the canonical bundle $K_{\Y_c}$ of the  Hermitian symmetric space of the compact type $\Y_c$ admits a $c_\Y$th-root. This is probably well-known to experts but we give some details because there are some subtleties we wish to discuss. A good general reference for what follows is~\cite{Kn}, see also~\cite{GS}. 

Group theoretically, $\Y_c=G_\C/Q$, where $Q$ is a maximal parabolic subgroup of the complexification $G_\C$ of $G$ such that $Q\cap G = K$, the maximal compact subgroup of $G$. If $U$ is a maximal compact subgroup of $G_\C$ containing $K$, then also $\Y_c=U/K$. Note that $G_\C$ is a connected simple complex Lie group. If $G_\C$ is not simply connected, let $\wt G_\C$ be its universal cover, and call $\wt Q$, $\wt U$, $\wt K$ the preimages of $Q$, $U$, $K$ in $\wt G_\C$. Then $\wt Q$, $\wt U$ and $\wt K$ are connected, $\wt Q$ is a maximal parabolic subgroup of $\wt G_\C$, $\wt U$ a maximal compact subgroup, and $\Y_c=\wt G_\C/\wt Q=\wt U/\wt K$.

\def\Aut{{\rm Aut}}
By a result of Murakami~\cite{Mu}, the Picard group $H^1(\Y_c,{\cal O}^\star)$ of $\Y_c$ is isomorphic to the character group $\Hom(\wt Q,\C^\star)$ of $\wt Q$. This can be seen as follows. First, every line bundle on $\Y_c$ is homogeneous, meaning that its automorphism group acts transitively on the base $\Y_c$. 
Indeed, the Picard group $H^1(\Y_c,{\cal O}^\star)$ of $\Y_c$ sits in the long exact sequence
$$
\cdots \fd H^1(\Y_c,{\cal O}) \fd H^1(\Y_c,{\cal O}^\star) \stackrel{c}{\fd} H^2(\Y_c,\Z) \fd H^2(\Y_c,{\cal O}) \fd\cdots
$$
where the connecting homomorphism $c$ assigns to the isomorphism class of a line bundle its first Chern class. By Kodaira vanishing theorem, because of the negativity of the canonical bundle of $\Y_c$, both $H^1(\Y_c,{\cal O})$ and $H^2(\Y_c,{\cal O})$ are reduced to zero, so $H^1(\Y_c,{\cal O}^\star)$ and $H^2(\Y_c,{\Z})$ are isomorphic. Now let $g$ be an element of $G_\C$ (or $\wt G_\C$) and let $L$ be any line bundle on $\Y_c$. Since $G_\C$ (or $\wt G_\C$) is connected and the Chern classes are integral classes, $c_1(g^\star L)=c_1(L)$ and hence $g^\star L$ and $L$ are isomorphic. 
Hence there is an automorphism of $L$ which acts by $g$ on $\Y_c$ and the transitivity is proved. In other words,  
the image $\bar G_\C$ of $\wt G_\C$ in $\Aut(\Y_c)$ is in the image of the natural morphism from $\Aut(L)$ to $\Aut(\Y_c)$. At the Lie algebra level, we hence have a surjective morphism ${\mathfrak a}{\mathfrak u}{\mathfrak t}(L)\fd\bar\gg$. Since $\bar\gg$ is simple, this map has a right inverse by Levi's theorem and, by simple connectedness, this gives a morphism from $\wt G_\C$ to $\Aut(L)$. Restricted to $\wt Q$, this morphism lands in the automorphism group of the fiber of $L$ above $e\wt Q\in\Y_c$, and thus defines a character $\chi$ of $\wt Q$. Moreover, the line bundle $L$ is just the bundle associated to the $\wt Q$-principal bundle $\wt G_\C\fd\Y_c$ via the action of $\wt Q$ on $\C$ by $\chi$.  

Denote by $\gg$, $\qg$, $\ug$, $\gg_0$ and $\kg_0$ the Lie algebras of $G_\C$, $Q$, $U$, $G$ and $K$. The complexification $\kg$ of $\kg_0$ has a unique $\ad\kg$-invariant complement $\pg$ in $\gg$. If we set $\pg_0=\gg_0\cap\pg$ then $\ug= \kg_0\oplus i\pg_0$ and $\gg_0=\kg_0\oplus\pg_0$ is a Cartan decomposition of $\gg_0$. Moreover, if $\hg_0$ is a maximal Abelian subalgebra of $\kg_0$, then it is a Cartan subalgebra of $\gg_0$ and its complexification $\hg$ is a Cartan subalgebra of $\gg$. This follows from the fact that, $\Y$ being Hermitian symmetric, $\kg_0$ has a 1-dimensional center $\zg_0$ which must be contained in $\hg_0$ and whose elements do not commute with non zero elements of $\pg_0$. Let $Z$ be the element of $\zg_0$ giving the complex structure of $\pg_0$ and hence of $\Y$. 
The complexification $\zg$ of $\zg_0$ satisfies $\zg\subset\hg\subset\kg\subset\qg$ and $\qg=\zg\oplus[\qg,\qg]$. 

Because of this last fact, the weight $d\chi_{|\hg}\in\hg^\star$ associated to a character $\chi$ of $Q$ or $\wt Q$ is entirely determined by its value on $Z$.  In particular, $\Hom(\wt Q,\C^\star)$ and hence $H^1(\Y_c,{\cal O}^\star)$ are isomorphic to $\Z$. 
So let $\L$ be the negative generator of $H^1(\Y_c,{\cal O}^\star)$ and consider the integer $n$ such that $K_{\Y_c}=\L^{\otimes n}$. We want to prove that $n=c_\Y$. Call $\lambda$ and $\kappa$ the weights in $\hg^\star$ corresponding to $\L$ and $K_{\Y_c}$, so that $\kappa=n\lambda$. 

Let $\Delta\subset\hg^\star$ be the set of roots of $\gg$ relative to the Cartan subalgebra $\hg$ and let $\Pi$ be a basis of $\Delta$. A root $\b\in\Delta$ is said to be compact if $\b(Z)=0$, noncompact otherwise. Let $\Delta_\kg$ be the set of compact roots, and $\Delta^+_\pg$ the set of positive noncompact roots. For each root $\b\in\Delta$, let $h_\b$ be the element of $\hg$ determined by $\b(Y)=B(h_\b,Y)$ for all $Y\in\hg$, where $B$ is the Killing form of $\gg$. Since $\ug$ is a compact real form of $\gg$, the roots are purely imaginary on $\hg_0$ and $h_\b\in i\hg_0$ for all $\b\in\Delta$.
The root space decomposition of $\gg$ is
$$
\gg=\hg\oplus\sum_{\b\in\Delta_\kg}\gg_\b\oplus\sum_{\b\in\Delta^+_\pg}\gg_{-\b}\oplus\sum_{\b\in\Delta^+_\pg}\gg_{\b}.
$$ 
Then $\kg = \hg\oplus\sum_{\b\in\Delta_\kg}\gg_\b$, $\pg^\pm:=\sum_{\b\in\Delta^+_\pg}\gg_{\pm\b}$ is the ($\pm i$)-eigenspace of $\ad(Z)$ on $\pg$ (in particular, $\b(Z)=\pm i$ if $\b\in\pm\Delta^+_\pg$), and $\qg=\kg\oplus\pg^-$.

Setting $d=\dim\Y_c=|\Delta^+_\pg|$, the canonical bundle of $\Y_c$ corresponds to the representation of $\wt Q$ on $\Lambda^d \pg^-$ coming from the representation of $\wt Q$ on $\pg^-$ induced by the adjoint representation of $\wt K$ on $\pg^-$.
By choosing a basis of $\pg^-$ adapted to the root space decomposition, we see that the action of $\ad(Y)$ on $\Lambda^d \pg^-$ for $Y\in\hg$ is simply given by $\kappa(Y)=-\sum_{\b\in\Delta^+_\pg}\b(Y)$. Therefore, $\kappa(Z)=-\sum_{\b\in\Delta^+_\pg}\b(Z)=-i\,d$. 

On the other hand, the weight $\lambda\in\hg^\star$ given by the generator $\L$ of $H^1(\Y_c,{\cal O}^\star)$ is minus the fundamental weight corresponding to the unique noncompact root $\zeta$ in the set of simple roots $\Pi$. More precisely, $\lambda$ is the weight defined by the relations $B(\lambda,\a)=\lambda(h_\a)=-\frac{1}{2}\|\a\|^2\delta_{\a,\zeta}$ for all $\a\in\Pi$ (here $\delta$ is the Kronecker delta and $\|\a\|^2:=B(\a,\a):=B(h_\a,h_\a)\geq 0$). Using the root space decomposition of $\gg$, it is clear that for all $Y\in\hg$, $B(Z,Y)=\tr\,(\ad(Z)\circ\ad(Y))=2i\sum_{\b\in\Delta^+_\pg}\b(Z)$ so that $Z=2i\sum_{\b\in\Delta^+_\pg}h_\b$. Therefore, 
$$
\lambda(Z)=2i\sum_{\b\in\Delta^+_\pg}\lambda(h_\b)=2i\sum_{\b\in\Delta^+_\pg}\lambda(h_{\zeta})
=-i\,d\,\|\zeta\|^2
$$
because if a root $\b=\sum_{\a\in\Pi} n_\a\a$ belongs to $\Delta^+_\pg$, then necessarily $n_{\zeta}=1$ since as we said ${\rm ad} (Z)$ acts by multiplication by $i$ on $\gg_\b$ for $\b\in\Delta^+_\pg$.
We conclude that $K_{\Y_c}=\L^{\otimes \|\zeta\|^{-2}}$. 

We claim that $\|\zeta\|^{-2}=c_\Y$. To prove this, we compare the $G$-invariant Kähler metric $\o_B$ on $\Y=G/K$ whose restriction to $T_{eK}\Y\simeq\pg_0$ is given by the Killing form $B$ and the initial Kähler metric $\o_\Y$ normalized so that its minimal holomorphic sectional curvature is $-1$. A classical computation shows that $c_1(K_\Y)=\frac{1}{4\pi}\o_B$, hence $\o_B=c_\Y\,\o_\Y$ and we are left with proving that the minimum of the holomorphic sectional curvature of $\o_B$ on $\Y$ is $-\|\zeta\|^2$.

We may choose a set $\Phi$ of $\rk\Y$ strongly orthogonal noncompact positive roots, meaning that for $\a,\b\in\Phi$, neither $\a+\b$ nor $\a-\b$ is a root, so that $B(\a,\b)=0$ for $\a\neq\b$.  Moreover we can find vectors $e_{\a}\in\gg_{\a}$ for $\a\in\pm\Phi$ such that $\bar e_\a= e_{-\a}$ (where complex conjugacy is w.r.t. the noncompact real form $\gg_0$ of $\gg$), $B(e_\a,e_\b)=\delta_{\a,-\b}$ and  $[e_\a,e_\b]= \delta_{\a,-\b}\,h_\a$. Since any vector of $\pg^+$ can be sent to an element in $\sum_{\a\in\Phi}\C e_\a$ by the adjoint action of $K$, it is enough to compute the holomorphic sectional curvature of a line spanned by such an element $Y=\sum_{\a\in\Phi}y_\a\,e_\a$, which is 
$$
\begin{array}{rcl}
\ds B\left(-\left[[Y,\bar Y],Y\right],\bar Y\right) & = & \ds -B\left([Y,\bar Y],[Y,\bar Y]\right)\\
& = & \ds - B\Big(\sum_{\a\in\Phi}|y_\a|^2\,h_\a,\sum_{\b\in\Phi}|y_\b|^2\,h_{\b}\Big)\\
& = & \ds - \sum_{\a,\b\in\Phi}|y_\a|^2|y_\b|^2B(h_\a,h_{\b})\\
& = & \ds - \sum_{\a\in\Phi}|y_\a|^4\,\|\a\|^2.
\end{array}
$$    
It follows from the classification of irreducible symmetric spaces that the root $\zeta$ coming from the complex structure on a Hermitian symmetric space is long (see~\cite[p. 414]{Kn}). Thus 
$$
B\left(-\left[[Y,\bar Y],Y\right],\bar Y\right) \geq - \|\zeta\|^2\,\|Y\|^4.
$$
Hence the minimum of the holomorphic sectional curvature of the Killing form is indeed $-\|\zeta\|^2$ and $c_\Y=\|\zeta\|^{-2}$ as claimed. Therefore $K_{\Y_c}=\L^{\otimes c_\Y}$. 

\medskip

Coming back to the definition of the Toledo invariant of the representation $\rho:\G\fd G$ and to the Milnor-Wood inequality, we denote by $\L_\Y$ the restriction to $\Y$ of the negative generator of $H^1(\Y_c,{\cal O}^\star)$. The canonical bundle $K_\Y$ of $\Y$ is simply the restriction to $\Y\subset\Y_c$ of the canonical bundle of $\Y_c$, so that $K_\Y=(\L_\Y)^{\otimes c_\Y}$. 
Given a $\rho$-equivariant map $f:\Bm\fd\Y$, we may pull-back $\L_\Y$ to a line bundle $f^\star\L_\Y$ on $\Bm$. Our discussion so far shows that, when the complexification $G_\C$ of $G$ is simply connected, $G$ acts by automorphisms on $\L_\Y$ and hence $\G$ acts by automorphisms on $f^\star\L_\Y$ which therefore descends to a line bundle on $X$. As before, the isomorphism class of this bundle is independent of the chosen equivariant map.   
The Toledo invariant of $\rho:\G\fd G$ is defined to be the degree of this bundle:

\begin{propdef}\label{defclassical}
Let $\G$ be a (torsion free) uniform lattice in the automorphism group of complex hyperbolic $m$-space $\Bm$. Let $\rho$ be a homomorphism of $\G$ into a linear connected simple noncompact Lie group of Hermitian type $G$ whose complexification is simply connected. Let $\Y$ be the Hermitian symmetric space of the noncompact type associated to $G$ and let $\Y_c$ be its compact dual. Then the pull-back of the restriction to $\Y$ of the negative generator of the Picard group of $\Y_c$ by any $\rho$-equivariant map $f:\Bm\fd\Y$ descends to a line bundle on $X=\G\backslash\Bm$ whose isomorphism class depends only on $\rho$ and which is denoted by $\L_\rho$. 

The Toledo invariant of $\rho$ is defined to be the $K_X$-degree of the bundle $\L_\rho$:   
$$\tau(\rho):=\deg \L_\rho=\L_\rho\cdot K_X^{m-1}.$$  
\end{propdef}

With this definition, $\tau(\rho)=\frac{1}{c_\Y}\deg f^\star K_\Y$ is a constant multiple of the Toledo invariant defined by Burger-Iozzi and their Milnor-Wood type inequality may be written:
$$
\left|\deg \L_\rho\right|\leq \rk G\,\frac{\deg K_X}{m+1}.
$$

\medskip

We end this section by discussing the additional condition on the group $G$ given in the definition, namely that its complexification $G_\C$ should be simply connected. In general, $G_\C$ is a finite cover of the group of automorphisms of the compact dual $\Y_c$ of the symmetric space associated to $G$. Here we have chosen the universal cover for this gives the nicest picture of the Toledo invariant as the degree of a line bundle. However, one often prefers to work with groups which admit faithful linear representations in dimension as small as possible. The classification of symmetric spaces for example is generally given in terms of those groups. Here is E. Cartan's list of irreducible Hermitian symmetric spaces $G/K$ of the noncompact type as found in~\cite{He}, where we inserted the value of the constant $c_{G/K}$ (for the two exceptional types, this value is taken from~\cite{Liu}):  

$$
\begin{array}{|c|c|c|c|c|c|c|}
\hline
\mbox{type} & G & K & G_\C & \rk G & \dim_\C G/K & c_{G/K} \\ 
\hline\hline
\mbox{A III} & \SU(p,q) & {\rm S}(\U(p)\times \U(q)) & \SL(p+q,\C) & \min(p,q) & pq & p+q\\ 
\hline
\mbox{BD I} & \SO_0(p,2) & \SO(p)\times\SO(2) & \SO(p+2,\C) & 2 & p & p\\ 
\hline
\mbox{C I} &{\rm Sp}(n,\R) & \U(n) & {\rm Sp}(n,\C) & n & n(n+1)/2 & n+1\\ 
\hline
\mbox{D III} & \SO^\star(2n) & \U(n) & \SO(2n,\C) & [n/2] & n(n-1)/2 & 2(n-1)\\ 
\hline
\mbox{E III} & & & & 2 & 16 & 12\\ 
\hline
\mbox{E VII} & & & & 3 & 27 & 18\\ 
\hline
\end{array}
$$

\bigskip

For our definition of the Toledo invariant, we asked for the line bundle $f^\star\L_\Y$ to descend on $X$. When $G_\C$ is not simply connected, this is not always possible because $G$ might not act on $\L_\Y$. More precisely, this is not true if some characters of the preimage $\wt Q$ of the maximal parabolic $Q<G_\C$ in $\wt G_\C$ do not come from characters of $Q$. The character groups of $Q$  and $K$ (resp. $\wt Q$ and $\wt K$) are isomorphic, so, in case $G_\C$ is not simply connected, it is enough to check whether all characters of $\wt K$ are lifts of characters of $K$.    
Going through the list of classical Hermitian symmetric spaces we just gave, we obtain: 
%\begin{itemize}
%\item 

\smallskip

-- type A III: the complexification of $G=\SU(p,q)$ is $\SL(p+q,\C)$ which is simply connected.
%\item 

\smallskip

-- type BD I: the complexification of $G=\SO_0(p,2)$ is $\SO(p+2,\C)$ whose universal cover is $\wt G_\C={\rm Spin}(p+2,\C)$. The maximal compact subgroup $K$ of $G$ is $\SO(p)\times\SO(2)$, its preimage in $\wt G_\C$ is $\wt K={\rm Spin}^c(p)={\rm Spin}(p)\times_{\Z_2}{\rm Spin(2)}$. If $\SO(2)$ and ${\rm Spin(2)}$ are identified with $\U(1)$, and if $\pi$ denotes the covering map ${\rm Spin}(p)\fd\SO(p)$, then the 2-sheeted covering map $\wt K\fd K$ is given by $(s,z)(=(-s,-z))\longmapsto(\pi(s),z^2)$. Since the characters of $\wt K$ are trivial when restricted to its derived subgroup ${\rm Spin}(p)\times\{1\}$, they all come from characters of $K$: $G$ acts on $\L_\Y$ and the bundle $f^\star\L_\Y$ goes down to $X$. 
%\item 

\smallskip

-- type C I: the complexification of $G={\rm Sp}(n,\R)$ is ${\rm Sp}(n,\C)$ which is simply connected.
%\item 

\smallskip

-- type D III: the complexification of $G=\SO^\star(2n)$ is $\SO(2n,\C)$ whose universal cover is $\wt G_\C={\rm Spin}(2n,\C)$. The maximal compact subgroup $K$ of $G$ is $\U(n)=\SU(n)\rtimes\U(1)$, its preimage $\wt K$ in $\wt G_\C$ is also isomorphic to a semidirect product $\SU(n)\rtimes\U(1)$ and the 2-sheeted covering map $\wt K\fd K$ is given by $(A,z)\longmapsto(A,z^2)$. In this case the triviality of the restrictions of the characters of $\wt K$ to its derived subgroup $\SU(n)\times\{1\}$ is not enough: if the character group of $\wt K$ is identified with $\Z$, only the characters belonging to $2\Z$ are lifts of characters of $K$. Therefore $G$ does not act on $\L_\Y$, only on $\L_\Y^{\otimes 2}$, so that only $f^\star\L_\Y^{\otimes 2}$ descends to a bundle $\L^2_\rho$ on $X$. Remark that $f^\star K_\Y$, seen as a line bundle on $X$, is nevertheless a power of $\L^2_\rho$, since in this case $c_\Y=2(n-1)$ is always even. The Toledo invariant could be defined as $\frac{1}{2}\deg_L\,\L^2_\rho$.             
%\end{itemize}

\section{The Toledo invariant: fundamental groups of varieties of general type}\label{general}

The reader will need some rudiments about the Minimal Model Program for which a good reference is \cite{KMM}. We only quote results which are useful for our purpose. 

\medskip

Let $X$ be a smooth projective variety of general type and of dimension $m\geq 2$. Recall that this means that the canonical bundle $K_X$ of $X$ is big, namely that $${\limsup}_{k\rightarrow\infty}\,\displaystyle\frac{\log({\rm dim}\,H^0(X,k K_X))}{\log k}=m~.$$

Let $\rho:\pi_1(X)\fd G\subset\GL(n,\C)$ be a representation of $\pi_1(X)$ in a connected simple noncompact Lie group of Hermitian type $G$ whose complexification is simply connected. Let $f:\widetilde X\fd\Y$ be a $\rho$-equivariant map from the universal cover of $X$ to the symmetric space associated to $G$, and $\L_\rho\fd X$ be (the isomorphism class of) the line bundle associated to $\rho$ we constructed in the previous section. 
As explained there, when $X$ is a ball quotient, the classical Toledo invariant of $\rho$ can be interpreted as the degree of $\L_\rho$, computed with respect to the polarization defined by the canonical class on $X$: $\tau(\rho)=\L_\rho\cdot(K_X)^{m-1}$. As a consequence, the Toledo invariant could be defined in the same way for representations of the fundamental group of any smooth variety with negative first Chern class. 

When the canonical bundle $K_X$ of $X$ is big but not nef, although we cannot consider the class $c_1(K_X)$ as a polarization, the number $\int_X  c_1(\L_\rho)\wedge c_1(K_X)^{m-1}$ makes sense and a first attempt would be to take it as a definition of the Toledo invariant.
The following simple example shows that this is too naive. Let $Y$ be a compact ball quotient of dimension 3, $C$ a smooth curve in $Y$, and let $f:X\fd Y$ be the blow-up of $Y$ along $C$. Let $\rho:\pi_1(X)\fd\SU(3,1)$ be the isomorphism $f_\star:\pi_1(X)\fd\pi_1(Y)$ induced by $f$. We want to compute the intersection number $\L_\rho\cdot K_X^2=\frac{1}{4}f^\star K_Y\cdot K_X^2$, where $K_Y$ is the canonical bundle of $Y$. We have $K_X=f^\star K_Y+E$, where $E$ is the exceptional divisor. Therefore,
$$
f^\star K_Y\cdot K_X^2 = (f^\star K_Y)^3 + f^\star K_Y\cdot E^2 + 2(f^\star K_Y)^2\cdot E = K_Y^3+f^\star K_Y\cdot E^2. 
$$    
Now, $f^\star K_Y\cdot E^2=(f^\star K_Y)_{|E}\cdot E_{|E}$ and, $E$ being isomorphic to the projectivization $\P(N_{C/Y})$ of the normal bundle $N_{C/Y}$ of $C$ in $Y$, $E_{|E}$ is the first Chern class of the tautological bundle ${\cal T}$ over $E$. Moreover, if ${K_Y}_{|C}=\sum n_i\,p_i$, then $(f^\star K_Y)_{|E}=\sum n_i\, E_{p_i}$ where $E_{p_i}$ is the fiber of $E\fd C$ over $p_i\in C$. Hence 
$$
f^\star K_Y\cdot E^2 = \sum n_i\int_{E_{p_i}}c_1({\cal T})_{|E_{p_i}}=(\sum n_i)\,\int_{\P^1} c_1({\cal O}(-1))=-K_Y\cdot C.
$$ 
If now we choose $C$ to be given by the zeros of two general sections of $K_Y^d$ for some $d$, we get  $K_Y\cdot C = d^2 K_Y^3$, so that $f^\star K_Y\cdot K_X^2 = (1-d^2)K_Y^3$ is unbounded in $d$ and we cannot hope for a Milnor-Wood inequality. 

\medskip

The generalization of the Toledo invariant to the case of smooth varieties of general type is therefore not straightforward. The previous example shows that we have to ``remove the negative part'' of $K_X$ in order to compute an invariant with good properties. In more precise terms, we need a Zariski decomposition of the canonical bundle of $X$, namely a decomposition $K_X=L+N$, where $L$ and $N$ are $\Q$-divisors, $L$ being nef and $N$ being effective, such that the natural map
$$H^0(X,{\cal O}_X(kK_X-\lceil kN\rceil))\fd H^0(X,{\cal O}_X(kK_X))$$
is bijective for any $k\geq 1$ (where $\lceil kN\rceil$ denotes the round up of the $\Q$-divisor $kN$).
In particular, this implies that we will have to work on special birational models of our given variety of general type, whose existence follows from the following deep result:

\begin{thm}[]{\rm (Birkar, Cascini, Hacon, McKernan~\cite{BCHM})}.
Let $X$ be a smooth variety of general type. Then, $X$ has a minimal model, which implies that the canonical ring $R(X,K_X)=\oplus_{k\in\N}H^0(X,{\cal O}(kK_X))$ is finitely generated. 
\end{thm}

Therefore the so-called canonical model $X_\ca:={\rm Proj}(R(X,K_X))$ of a smooth variety of general type $X$ is well defined. The variety $X$ is birational to its canonical model $X_\ca$ but there is maybe no morphism $X\fd X_\ca$. However, if $\varphi: X'\fd X_\ca$ is a resolution of singularities of $X_\ca$, then $X'$ is a birational model of $X$ of the type we are looking for. Indeed, the canonical model $X_\ca$ has canonical singularities, meaning that it is a normal variety, its canonical divisor $K_{X_\ca}$ is a  $\Q$-Cartier divisor, and if $\varphi: X'\fd X_\ca$ is as above, $K_{X'}=\varphi^\star K_{X_\ca}+N$ with $N$ effective and supported in the exceptional divisor. Moreover, the canonical divisor $K_{X_\ca}$ is ample, so that $\varphi^\star K_{X_\ca}$ is indeed nef. 
Conversely, by the base-point-free theorem in its generalized form~\cite{Ka}, if $K_{X'}$ admits a Zariski decomposition $K_{X'}=L+N$ for some birational model $X'$ of $X$, then $L$ is semi-ample which implies that there exists a proper modification $\varphi:X'\fd X_\ca$ and that $L=\varphi^\star K_{X_\ca}$ as above. We call such a smooth model of $X$ a {\it good model}. Even if varieties of general type are not necessarily projective, a good model is always projective.

\medskip

As is well known, two smooth birational varieties have isomorphic fundamental groups. Therefore, if $\rho$ is a representation of the fundamental group of a smooth variety of general type and if $X$ is a good model of this variety, we can consider $\rho$ as a representation of $\pi_1(X)$. This allows us to give the following

\begin{defi}\label{defgeneral}
Let $\G$ be the fundamental group of a smooth variety of general type and dimension $m\geq 2$ and let $X_\ca$ and $X\stackrel{\varphi}{\fd}X_\ca$ be respectively the canonical model and a good model of this variety.
Let $\rho$ be a representation of $\G$ in a linear connected simple noncompact Lie group of Hermitian type $G$ whose complexification is simply connected. 

The Toledo invariant of $\rho$ is the $(\varphi^\star K_{X_\ca})$-degree of the line bundle $\L_\rho\fd X$ associated to $\rho$ as in Section~\ref{classical}:
$$
\tau(\rho)=\deg_{\varphi^\star K_{X_\ca}}(\L_\rho)=\L_\rho\cdot(\varphi^\star K_{X_\ca})^{m-1}=\int_{X}c_1(\L_\rho)\wedge c_1(\varphi^\star K_{X_\ca})^{m-1}.
$$
\end{defi}

It is easy to see that the invariant is well-defined, namely that it does not depend on the choice of the good model. Indeed, let $X$ and $X'$ be two good models of the same variety of general type. Then, by a theorem of Hironaka, there exists a manifold $X_0$ and proper modifications $\nu:X_0\fd X$ and $\nu':X_0\fd X'$  ($X_0$ is automatically a good model),
such that the following diagram commutes:
$$\xymatrix{&X_0\ar[dl]_{\nu}\ar[dr]^{\nu'} &\\
X\ar[dr]_\varphi\ar@{-->}[rr]&& X'\ar[dl]^{\varphi'}\\
& X_\ca&}$$
If $f:\widetilde{X}\fd \Y$ and $f':\widetilde{X'}\fd \Y$ are smooth $\rho$-equivariant maps (keeping in mind that the fundamental groups of $X$, $X'$ and $X_0$ are naturally isomorphic) then $f\circ\tilde\nu$ and $f'\circ\tilde\nu'$ are $\rho$-equivariant maps on $\widetilde{X_0}$ (where $\tilde\nu$ and $\tilde\nu'$ denote lifts of $\nu$ and $\nu'$ on the respective universal covers). Then, it is immediate that
$$
\begin{array}{rcl}
c_1(\L_\rho)\cdot c_1(\varphi^\star K_{X_\ca})^{m-1}&=&c_1(\nu^\star \L_\rho)\cdot c_1(\nu^\star\varphi^\star K_{X_\ca})^{m-1}\\
& = & c_1({\nu'}^\star \L_\rho')\cdot c_1({\nu'}^\star{\varphi'}^\star K_{X_\ca})^{m-1}\\
& = & c_1(\L_\rho')\cdot c_1({\varphi'}^\star K_{X_\ca})^{m-1}~.
\end{array}
$$

\medskip

We remark also that ${\rm deg}_{\varphi^\star K_{X_\ca}}\,(K_{X})=K_{X}\cdot\varphi^\star K_{X_\ca}^{m-1}$ is independent of the good model $X$ of the variety we are considering and that its value is $K_{X_\ca}^m$.

\section{Proof of the main result}\label{proof}

We recall the statement of our main theorem:

\begin{theo}\label{mainthm}
Let $X$ be a smooth variety of general type and of complex dimension $m\geq 2$, and let $X_\ca$ be its canonical model.  
Let $G$ be either $\SU(p,q)$ with $1\leq q\leq 2\leq p$, ${\rm Spin}(p,2)$ with $p\geq 3$, or ${\rm Sp}(2,\R)$. 
Finally let $\rho:\pi_1(X)\fd G$ be a representation. 

Then the  Toledo invariant of $\rho$ satisfies the Milnor-Wood type inequality  
$$
|\tau(\rho)|\leq \rk G\,\frac{K_{X_\ca}^m}{m+1},
$$
where $K_{X_\ca}$ is the canonical divisor of $X_\ca$.  

Equality holds if and only if $G=\SU(p,q)$ with $p\geq qm$ and there exists a $\rho$-equivariant (anti)holomorphic proper embedding from the universal cover of $X_\ca$ onto a totally geodesic copy of complex hyperbolic $m$-space $\Bm$, of induced holomorphic sectional curvature $-1/q$, in the symmetric space associated to $G$. In particular, $X_\ca$ is then smooth and uniformized by $\Bm$, and the representation $\rho$ is discrete and faithful.
\end{theo}

We will first prove the inequality, and then study the equality case. In what follows, the smooth variety of general type $X$ that appears in the theorem will be assumed to be a good model, i.e. to be such that the pluricanonical map $\varphi:X\fd X_\ca$ is a morphism. As explained before, there is no loss of generality in doing so. Unless otherwise specified, all degrees on $X$ will be computed with respect to the polarization $\varphi^\star K_{X_\ca}$ and we will abbreviate $\deg_{\varphi^\star K_{X_\ca}}$ to $\deg$.

\subsection{Proof of the Milnor-Wood type inequality}\label{inequality}

We recall briefly the notions of semistability that we will need: let $L$ be a nef $\Q$-line bundle on a complex manifold $X$ and let ${\cal F}$ be a coherent (saturated) sheaf of positive rank on $X$. We already defined the $L$-degree ${\rm deg}_L {\cal F}=c_1({\cal F})\cdot c_1(L)^{m-1}$ of ${\cal F}$. Its $L$-slope is defined by $\mu_L({\cal F})=\frac{{\rm deg}_L{\cal F}}{\rk {\cal F}}$.
We have the corresponding notion of $L$-semistability of a holomorphic vector bundle $E$ on $X$: this means that for any coherent subsheaf ${\cal F}$ of ${\cal O}_X(E)$, $\mu_L({\cal F})\leq\mu_L(E)$ if $\rk({\cal F})>0$.

Similarly, a Higgs bundle $(E,\theta)$ on $X$ (that is a holomorphic vector bundle $E$ together with a holomorphic map $\theta:E\fd E\otimes\Omega^1_X$ such that $\theta\wedge\theta=0$) is called $L$-semistable if the inequality $\mu_L({\cal F})\leq\mu_L(E)$ holds for any coherent Higgs subsheaf of positive rank ${\cal F}$ of ${\cal O}_X(E)$, namely for any subsheaf of positive rank ${\cal F}$ of ${\cal O}_X(E)$ such that $\theta({\cal F})\subset {\cal F}\otimes \Omega^1_X$.
The Higgs bundle is said to be $L$-stable if the inequality is strict when $\rk({\cal F})<\rk(E)$. Finally, we say that $(E,\theta)$ is $L$-polystable if it is a direct sum of stable Higgs bundles of the same slope.

\medskip

The following result will be crucial for the proof of Theorem~\ref{mainthm}:

\begin{thm}[]{\rm (Enoki~\cite{En})}.
Let $X\stackrel{\varphi}{\fd} X_\ca$ be a good model of a variety of general type. Then the tangent bundle $T_{X}$ of $X$ is $(\varphi^\star K_{X_\ca})$-semistable.
\end{thm}

Let now $\rho:\pi_1(X)\fd\GL(n,\C)$ be a reductive representation. By a result of Corlette~\cite{C}, there exists a harmonic $\rho$-equivariant map $f:\widetilde X\fd \Y$, where $\Y=\GL(n,\C)/\U(n)$. This map allows to construct a Higgs bundle $E$ of rank $n$ and degree 0 on $X$ (through the standard representation of $\GL(n,\C)$ on $\C^n$) which is polystable with respect to any polarization coming from an ample divisor, see~\cite{S1}. Although the divisor $\varphi^\star K_{X_\ca}$ is not necessarily ample, we have the

\begin{lemma}\label{polystab}
Let $X\stackrel{\varphi}{\fd} X_\ca$ be a good model of a smooth variety of general type. Let $\rho:\pi_1(X)\fd \GL(n,\C)$ be a reductive representation. Then, the associated Higgs bundle $(E,\theta)$ is $(\varphi^\star K_{X_\ca})$-polystable.
\end{lemma}

\begin{demo}
As in \cite{S1}, this is a consequence of the Chern-Weil formula, together with the fact that, thanks to the ampleness of $K_{X_\ca}$, $\varphi^\star K_{X_\ca}$ can be endowed with a smooth metric whose curvature is semi-positive, and strictly positive outside the exceptional set ${\rm Ex}(\varphi)$ of $\varphi$. We give details of the proof of the lemma for the sake of completeness.

The flat connection $D_E$ on $E$ can be written $D_E=D_H+\theta+\theta^\star$ where $D_H$ is the component preserving the harmonic metric $H$ defined by $f$.
Let ${\cal F}\subset{\cal O}_X(E)$ be a saturated Higgs subsheaf and let $Y\subset X$ be an analytic subset of $X$ of codimension at least 2 such that ${\cal F}$ is a vector subbundle of $E$ outside of $Y$. On $X\backslash Y$, we can decompose the flat connection $D_E$ with respect to the orthogonal decomposition $E={\cal F}\oplus{\cal F}^\perp$ where the background metric is $H$, and denoting by $\sigma\in C^\infty_{1,0}(X,{\rm Hom}({\cal F},{\cal F}^\perp))$ the second fundamental form of ${\cal F}$:
$$D_E=\left(\begin{matrix}D_{\cal F} & -\sigma^\star\\
\sigma & D_{{\cal F}^\perp}\\
  \end{matrix}\right)+\left(\begin{matrix} \theta_1+\theta_1^\star & \theta_2\\
  \theta_2^\star & \theta_3+\theta_3^\star\\
 \end{matrix}\right)~.
$$
Since $D_E^2=0$, we have $(D_{\cal F}+\theta_1+\theta_1^\star)^2=-(\theta_2-\sigma^\star)\wedge(\theta_2^\star+\sigma)$. In order to compute the degree of ${\cal F}$, we can use the connection $D_{\cal F}+\theta_1+\theta_1^\star$ and we obtain
$$
\deg{\cal F}=\frac{i}{2\pi}\int_{X\backslash Y}{\rm tr}\,(-\theta_2\wedge\theta_2^\star+\sigma^\star\wedge \sigma)\wedge c_1(\varphi^\star K_{X_\ca})^{m-1}\leq 0~.
$$
If moreover $\deg {\cal F}=0$, then $\theta_2$ and $\sigma$ vanish identically on $X\backslash (Y\cup {\rm Ex}(\varphi))$ (because $c_1(\varphi^\star K_{X_\ca})$ is strictly positive on $X\backslash {\rm Ex}(\varphi)$) and hence they vanish on $X\backslash Y$ by continuity. Therefore $E$ splits holomorphically as a direct sum of two Higgs subbundles on $X\backslash Y$. A classical argument due to Lübke (see~\cite{K}, p.~179) shows that this splitting extends to $X$.
\end{demo}

\medskip

When the representation $\rho$ takes values in a Lie group of Hermitian type, the associated Higgs bundle $E$ has some extra-structure, see for example~\cite{KM}. Since representations $\rho:\pi_1(X)\fd G=\SU(p,q)\subset \SL(p+q,\C)$, $p\geq q\geq 1$ will be of particular importance for us, we give some more details in this case. The symmetric space $\Y$ associated to $G=\SU(p,q)$ has rank $q$ and is an open subset in the Grassmannian $\Y_c$ of $q$-planes in $\E=\C^{p+q}$. If $\W\in\Y$ is a $q$-plane and if $\V=\E/\W$, the holomorphic tangent bundle of $\Y$ at $\W$ is identified with $\Hom(\W,\V)$, which we see as the space of $q\times p$ complex matrices. With respect to this identification, the Hermitian product $h$ coming from the Kähler metric $\o_\Y$ is given by $h(A,B)=4\,\tr\trans\bar B A$, so that the holomorphic sectional curvature of the complex line spanned by $A$ equals $-\tr(\trans\bar AA)^2/(\tr\trans\bar AA)^2$, which is indeed pinched between $-1$ and $-\frac{1}{q}$. The trivial bundle $\underline{E}=\Y_c\times\E$ over $\Y_c$ and its tautological subbundle $S$ restrict to vector bundles on $\Y$, denoted by the same letters. 
The bundle $S$, resp. the quotient bundle $\underline{E}/S$, pulls back to a rank $q$, resp. $p$, holomorphic bundle $W$, resp. $V$, on $X$. The Higgs bundle $E$, which as a smooth bundle is the flat bundle obtained by pulling back $\underline{E}$ to $X$, splits holomorphically as the sum $V\oplus W$. The Higgs field $\theta\in \Omega^1(X,{\rm End}\, E)$ is in fact a holomorphic 1-form with values in ${\rm Hom}\,(W,V)\oplus{\rm Hom}\,(V,W)$. We shall denote by $\beta$ (resp. $\gamma$) the projection of $\theta$ on the first (resp. second) summand of this decomposition. The decomposition $V\oplus W$ is orthogonal for the harmonic metric of $E$. The curvature forms of $V$ and $W$ w.r.t. this metric are given by
$$
R^V = -\b\wedge\b^\star - \g^\star\wedge\g\quad\mbox{ and }\quad
R^W = -\b^\star\wedge\b - \g\wedge\g^\star.
$$
Moreover, the bundle $f^\star T_\Y$ which goes down on $X$ is isomorphic to $\Hom(W,V)$ so that $\deg(f^\star K_\Y)=-\deg({\rm Hom}\,(W,V))=p\,\deg W-q\,\deg V=(p+q)\,\deg W$. The line bundle $\L_\rho$ is the pull-back of the determinant bundle of the tautological bundle on the Grassmannian, i.e. $\L_\rho=\det W$.  
 
\medskip

We are now ready to prove the 

\begin{prop}\label{ineq}
Let $X$ be a smooth projective variety of general type of dimension $m$ and let $X_\ca$ be its canonical model. Let $G$ be either $\SU(p,q)$ with $1\leq q\leq 2\leq p$, ${\rm Spin}(p,2)$ with $p\geq 3$ or ${\rm Sp}(2,\R)$. Finally let $\rho:\pi_1(X)\fd G$ be a representation. Then the Milnor-Wood type inequality
$$|\tau(\rho)|\leq \rk G\,\frac{K_{X_\ca}^m}{m+1}$$
holds.
\end{prop}

\begin{demo}
In fact the proof is exactly the same as in the classical case where $X$ is a compact quotient of $\Bm$, the key ingredients being only the semistability of the Higgs bundle $E$ associated to $\rho$ (Lemma~\ref{polystab}) and the classical semistability of $T_X$ (Enoki's theorem) with respect to the  polarization chosen to define the Toledo invariant. Since the proof is much more involved for higher rank Lie groups, where one needs to make use of the $\C^\star$-action on the moduli space of polystable Higgs bundles, we refer the reader to~\cite{KM} for rank two $G$'s and here we illustrate the importance of the stability properties by treating the case $G=\SU(p,1)$.

%We only prove the inequality for $G=\SU(p,1)$, whose associated symmetric space $\Y$ is complex hyperbolic $p$-space. When $G$ is one of the other groups listed in %Proposition~\ref{ineq}, the proof given in~\cite{KM} (for $X$ a compact quotient of $\Bm$) works equally well to yield the inequality, the key ingredients being %only the semistability of $E$ and $T_X$ with respect to the chosen polarization. Compared to the case of $\SU(p,1)$, the additional difficulty comes from the need %of making use of the $\C^\star$-action on the moduli space of polystable Higgs bundles.

Assume first that $\rho$ is not reductive. Then there exists a 1-form $\varsigma$ on $\Y=\B^p$ which is invariant by the action of $\rho(\pi_1(X))$ and such that $\d\varsigma=\frac{1}{4\pi}\,\o_\Y$, see for example~\cite{KM0}. Therefore, $\tau(\rho)=\int_{X} \df^\star\varsigma\wedge c_1(\varphi^\star K_{X_\ca})^{m-1}=0$.

We suppose now that $\rho:\pi_1(X)\fd \SU(p,1)$ is reductive and we call $E=V\oplus W$ the Higgs bundle associated to $\rho$ defined above.
%The bound on $\deg(f^\star K_Y)$ easily follows from the semistability of $E$ as a Higgs bundle (Lemma~\ref{polystab}), and the classical semistability of $T_X$ %(Enoki's theorem).
Since $T_{X}$ is a semistable vector bundle, its twist by the line bundle $W$ is semistable too and the image $\Im\b$ of the morphism $\beta:W\otimes T_X\fd V$ thus satisfies $\mu(W)+\mu(T_{X})\leq\mu(\Im\b)$ (remark that if $\beta\equiv 0$, $W$ is a Higgs subsheaf of $E$ and we immediately get $\deg W\leq0$). Combining this inequality with the fact that $W\oplus \Im\b$ is a Higgs subsheaf of $E$, we obtain:
$$\deg W+\frac{\deg T_{X}}{m}\leq \frac{\deg \Im\b}{r}\leq-\frac{\deg W}{r} 
$$
where $r$ is the generic rank of $\beta$. Therefore, we get the desired bound
$$\deg W\leq\frac{r}{r+1}\frac{\deg\Omega^1_X}{m}\leq \frac{\deg\Omega^1_X}{m+1}=\frac{K^m_{X_\ca}}{m+1} 
$$
because $r\leq m$.

The proof of the lower bound goes exactly the same way by considering the dual of the Higgs bundle $E$.
\end{demo}

\medskip

\subsection{Maximal representations}\label{equality}

In this section all representations will be assumed to have a positive Toledo invariant. This is not a loss of generality since if a representation is maximal with negative Toledo invariant, we may change the complex structure of $X$ for its opposite to get our result. 

\medskip

If we look back at the proof of the Milnor-Wood inequality we just gave in the case $G=\SU(p,1)$, we see that if $\tau(\rho)$ is maximal, then $\b:W\otimes T_X\fd V$ is generically injective and $\deg W+\deg \Im\b=0$. This will be a crucial information and it follows from \cite{KM} that the same kind of conclusion holds for maximal representations in the Lie groups listed in Theorem~\ref{mainthm}. To be more precise, if $\rho$ is a maximal representation of the fundamental group of a variety of general type and dimension $m\geq 2$ in such a Lie group $G$, then necessarily $\rho$ is reductive, $G$ is $\SU(p,q)$ with $q=1,2$ and $p\geq qm$, and moreover if $(E=V\oplus W,\theta=\b\oplus\g)$ is the Higgs bundle associated to $\rho$ on a good model $X$, then $\deg W+\deg \Im\b=0$ and $\b:W\otimes T_X\fd V$ is generically injective. 

We believe that this is still true in the general case where no restriction is made on the rank of the target Lie group $G$. Therefore, since the rest of our arguments are valid in this more general setting, we formulate the following proposition, which implies our main theorem.

\begin{prop}
Let $\rho$ be a reductive representation of the fundamental group of a variety of general type and dimension $m\geq 2$ in the group $G=\SU(p,q)$ with $p\geq qm$ and $q\geq 1$. Let $X_\ca$ and $X\stackrel{\varphi}{\fd}X_\ca$ be respectively the canonical model and a good model of this variety, and let $E=V\oplus W\fd X$ be the $(\varphi^\star K_{X_\ca})$-polystable $G$-Higgs bundle associated to the $\rho$-equivariant harmonic map $f$ from the universal cover $\widetilde X$ of $X$ to the Hermitian symmetric space $\Y$ associated to $G$. 

Assume that $\rho$ is maximal, i.e. $\deg W=q\frac{K^m_{X_\ca}}{m+1}$, and assume moreover that $\deg W+\deg \Im\b=0$ and that $\b: W\otimes T_X\fd V$ is generically injective. 

Then $f:\widetilde X\fd\Y$ factors through a $\rho$-equivariant holomorphic proper embedding of the universal cover of $X_\ca$ into a totally geodesic copy of complex hyperbolic $m$-space $\Bm$ in $\Y$, of induced holomorphic sectional curvature $-\frac{1}{q}$. In particular, the canonical model $X_\ca$ is smooth and uniformized by $\Bm$, and the representation $\rho$ is faithful and discrete.
\end{prop}

\begin{demo}
The fact that $\b: W\otimes T_X\fd V$ is generically injective implies of course that the morphism $\df:T_X\fd W^\star\otimes V\simeq f^\star T_\Y$ is generically injective, namely that the $\rho$-equivariant harmonic map $f:\widetilde X\fd\Y$ is a generic immersion. But it has much stronger consequences. 

The first one is that $\gamma\equiv 0$, which is equivalent to the holomorphicity of $f$. We give a short proof of that point: let $x\in X$ be such that $\beta_x$ is injective and let $\xi,\eta\in T_x^1$ be two linearly independent tangent vectors (recall that $m\geq 2$). The relation $\theta\wedge\theta=0$ implies that $\b_x(\xi)\gamma_x(\eta)v=\b_x(\eta)\gamma_x(\xi)v$ for any $v\in V_x$. Since $\b_x$ is injective, we must have $(\gamma_x(\eta)v)\otimes \xi=(\gamma_x(\xi)v)\otimes \eta$ hence $\gamma_x(\eta)v=\gamma_x(\xi)v=0$ and this is true for any $\xi,\eta$ and $v$ as above, thus $\gamma_x\equiv 0$. As a consequence, $\gamma$ vanishes generically on $X$ hence identically.

At this point it is already possible to deduce that $f:\widetilde X\fd\Y$ factors through the universal cover $\widetilde X_\ca$ of the canonical model of $X$. This follows from the  

\begin{thm}[]{\rm (Takayama~\cite{Ta})}. Let $X$ be a smooth variety of general type and $X_\ca$ its canonical model. Then the fundamental groups $\pi_1(X)$ and $\pi_1(X_\ca)$ are naturally isomorphic (meaning that the isomorphism is induced by the map between the varieties).
\end{thm}

Indeed, let $A\subset X_\ca$ be a subvariety of codimension $\geq 2$ such that $X_\ca\backslash A$ is biholomorphic to $X\backslash B$ where $B=\varphi^{-1}(A)$. Let us denote by $\widetilde A$ resp. $\widetilde B$ the lift of $A$ resp. $B$ in $\widetilde X_\ca$ resp. $\widetilde X$ and $\tilde\varphi$ a lift of $\varphi$ to the universal covers. Since $\pi_1(X)$ and $\pi_1(X_\ca)$ are isomorphic, the restriction $\tilde \varphi_{|\widetilde X\backslash \widetilde B}:\widetilde X\backslash \widetilde B\fd \widetilde X_\ca\backslash\widetilde A$ is a biholomorphism too. Finally, the holomorphic map $f\circ(\tilde\varphi_{|\widetilde X\backslash \widetilde B})^{-1}:\widetilde X_\ca\backslash\widetilde A\fd\Y\subset\C^{pq}$ can be extended to a $\rho$-equivariant holomorphic map $g:\widetilde X_\ca\fd\Y$ by normality of $\widetilde X_\ca$, and $f=g\circ\tilde\varphi$. 

\medskip

The degree of the Higgs subsheaf $W\oplus{\rm Im}\b$ of $E$ being zero, if we call $V'$ the saturation of ${\rm Im}\b$ in $V$, we get by polystability of $E$ that $V'$ is a holomorphic subbundle of $V$ of rank $qm$ and that there exists a holomorphic subbundle $V''$ of $V$ such that $E$ splits as the sum $(W\oplus V')\oplus V''$ of two Higgs subbundles. The Higgs field of $E$ is reduced to $\b$, which we now see as a morphism from $W\otimes T_X$ to $V'$. From the curvature formulas given in Subsection~\ref{inequality}, we have $R^{V'}=-\b\wedge\b^\star$ and $R^{W}=-\b^\star\wedge\b$, whereas $V''$ is a flat unitary bundle. In particular, $c_1(V') = -c_1(W) = {\rm tr}(\b^\star\wedge\b) =-\frac{1}{4\pi}f^\star\o_\Y$.

The Higgs field $\b:W\otimes T_X\fd V'$ is generically injective, and so is $\det\b:(\det W)^{m}\otimes K_X^{-q}\fd \det V'$. Hence there exists an effective divisor $D$ on $X$ such that the line bundles $(\det W)^{m}\otimes K_X^{-q}\otimes[D]$ and $\det V'$ are isomorphic outside codimension 2. The maximality of the Toledo invariant says that the line bundles $(\det W)^{m}\otimes K_X^{-q}$ and $\det V'$ have the same degree, so that necessarily $\deg [D]=0$, i.e. $\int_D c_1(\varphi^\star K_{X_\ca})^{m-1}=0$. Since $K_{X_\ca}$ is ample, the support of $D$ is included in the exceptional set ${\rm Ex}(\varphi)$ of $\varphi$. In particular, $f:\widetilde X\fd\Y$ is an immersion outside of the union of the lift of ${\rm Ex}(\varphi)$ to $\wt X$ with a codimension 2 subset of $\wt X$.

\medskip

We now want to prove that the image of $f$ lies in a totally geodesic copy of the ball $\Bm$ in the symmetric space $\Y$ associated to $\SU(p,q)$.
As indicated to us by Mok, if $\rho$ is discrete and $f$ is an embedding, this follows from Theorem 3 in his joint paper \cite{EM} with Eyssidieux. Here we extend their arguments, although in a somewhat different manner, to the generality we need. The key idea is to use again the injectivity of $\b$, this time in computing degrees with respect to the polarization $f^\star\o_\Y$. 

Consider the morphism $\df:T_X\fd W^\star\otimes V'$ and call $I$ the saturation of its image in $W^\star\otimes V'$. The sheaf $I$ is a subbundle of $W^\star\otimes V'$ outside codimension 2 and $\df:T_X\fd I$ is generically an isomorphism. As before, there exists an effective divisor $D'$ such that $\det\df:K_X^{-1}\otimes [D']\fd \det I$ is an isomorphism outside codimension 2. This implies that
$$
c_1(I)\cdot(f^\star\o_\Y)^{m-1} = c_1(K_X^{-1}\otimes [D'])\cdot(f^\star\o_\Y)^{m-1}.
$$
Note that the support of $D'$ is included in the support of $D$, hence in ${\rm Ex}(\varphi)$, so that
$$
c_1(K_X^{-1}\otimes [D'])\cdot(f^\star\o_\Y)^{m-1}= c_1(K_X^{-1})\cdot(f^\star\o_\Y)^{m-1}=\frac{1}{q}\,c_1(K_X^{-q}\otimes [D])\cdot(f^\star\o_\Y)^{m-1},
$$
for $f=g\circ\tilde\varphi$. Now,
$$
c_1(K_X^{-q}\otimes [D])\cdot(f^\star\o_\Y)^{m-1}=c_1(\det V'\otimes(\det W)^{-m})\cdot(f^\star\o_\Y)^{m-1}=-\frac{m+1}{4\pi}\int_X(f^\star\o_\Y)^{m}
$$
because  $K_X^{-q}\otimes[D]$ and $\det V'\otimes(\det W)^{-m}$ are isomorphic outside codimension 2.
Hence
$$
c_1(I)\cdot(f^\star\o_\Y)^{m-1} = -\frac{1}{2\pi}\,\frac{m+1}{2q}\int_X(f^\star\o_\Y)^{m}.
$$

On the other hand, $c_1(I)\cdot(f^\star\o_\Y)^{m-1}$ can be computed using the fact that $I$ is a subbundle of $W^\star\otimes V'$ (outside codimension 2). If $\sigma$ is its second fundamental form in $W^\star\otimes V'$, then the first Chern form of $I$ is
$$
\frac{i}{2\pi}\tr (R^I) = \frac{i}{2\pi}\Big(\tr(R^{W^\star\otimes V'}_{|I})-\tr(\sigma\wedge \sigma^\star)\Big).
$$
Around a point $x$ in the dense open subset of $X$ where $\df:T_X\fd I\subset W^\star\otimes V'$ is injective we choose local coordinates $\{z_k\}$ such that $\b_k:=\df(\frac{\partial}{\partial z_k})$ is an orthonormal basis of the fiber of $I$ at $x$ with respect to the metric $f^\star \o_\Y$. We identify the fiber of $W^\star\otimes V'$ with the space of $q\times qm$ complex matrices endowed with the pull-back of the Hermitian scalar product $h$ on $T_\Y$, see Subsection~\ref{inequality}. 
Then, $f^\star\o_\Y=\frac{i}{2}\sum_k dz_k\wedge d\bar z_k$ and
$R^{W^\star\otimes V'}\xi=-\sum_{j,k}(\b_j\trans\bar\b_k\xi+\xi\trans\bar\b_k \b_j)dz_j\wedge d\bar z_k$, so that
$$
\tr(R^{W^\star\otimes V'}_{|I})=-4\sum_{j,k,l}\tr(\trans\bar\b_l\b_j\trans\bar\b_k\b_l+\trans\bar\b_l\b_l\trans\bar\b_k \b_j)dz_j\wedge d\bar z_k.
$$
Hence,
$$
\begin{array}{rcl}
\ds i\tr(R^{W^\star\otimes V'}_{|I})\wedge (f^\star\o_\Y)^{m-1} & = &
\ds -\frac{8}{m} \Big(\sum_{k,l}\tr(\trans\bar\b_l\b_k\trans\bar\b_k\b_l+\trans\bar\b_l\b_l\trans\bar\b_k \b_k)\Big)(f^\star\o_\Y)^{m}\\
& \leq &
\ds -\frac{8}{m} \Big(\sum_{k}\tr(\trans\bar\b_k\b_k)^2+\tr\big(\sum_k\trans\bar\b_k\b_k\big)^2\Big)(f^\star\o_\Y)^{m}\\
& \leq &
\ds -\frac{8}{qm} \Big(\sum_{k}(\tr\trans\bar\b_k\b_k)^2+\big(\sum_k\tr\trans\bar\b_k\b_k\big)^2\Big)(f^\star\o_\Y)^{m}\\
%& = &
%\ds -\frac{2}{qm}\frac{m(m+1)}{4}(f^\star\o_\Y)^{m}\\
& = &
\ds -\frac{m+1}{2q}(f^\star\o_\Y)^{m}.
\end{array}
$$
By continuity, this inequality holds at the points where both sides are defined, namely outside a codimension 2 subset of $X$.

Summing up, we have
$$
\begin{array}{rcl}
\ds -\frac{m+1}{2q}\int_X(f^\star\o_\Y)^{m} & = & \ds (2\pi)\,c_1(I)\cdot(f^\star\o_\Y)^{m-1} \\
& \leq & \ds -\frac{m+1}{2q}\int_X(f^\star\o_\Y)^{m}-\int_X i\tr(\sigma\wedge \sigma^\star)\wedge (f^\star\o_\Y)^{m-1}.
\end{array}
$$
Therefore the semipositive form $i\tr(\sigma\wedge \sigma^\star)\wedge (f^\star\o_\Y)^{m-1}$ vanishes, and equality holds in the inequality
$i\tr(R^{W^\star\otimes V'}_{|I})\wedge (f^\star\o_\Y)^{m-1}\leq -\frac{m+1}{2q}(f^\star\o_\Y)^{m}$. The first point means that the second fundamental form of $I$ is zero on the dense open set where it is defined and $f^\star\o_\Y$ is positive definite, so that $f$ maps $\widetilde X$ into a totally geodesic $m$-dimensional submanifold ${\mathcal Z}$ of $\Y$. The second one implies that for all $k\neq l$, all the column vectors of $\b_k$ are orthogonal to all the column vectors of $\b_l$, and that for all $k$ the column vectors of $\b_k$ are pairwise orthogonal and have the same norm (w.r.t. the standard Hermitian scalar product on $\C^{qm}$). Using the formula for the holomorphic sectional curvature in $\Y$, we see that the holomorphic sectional curvature of every complex line in $\df(T_{\widetilde X})\subset T_\Y$ equals $-\frac{1}{q}$. Hence the totally geodesic submanifold ${\mathcal Z}$ is indeed a ball $\Bm$ of maximal possible holomorphic sectional curvature in $\Y$.

\medskip

We are ready to prove that $X_\ca$ is smooth. 
We may consider that the target of the maps $f$ and $g$ is $\mathcal Z$, and that they are equivariant with respect to the induced representation of $\pi_1(X)$ in the automorphism group ${\rm Aut}(\mathcal Z)$ of $\mathcal Z$, which will be still denoted by $\rho$. Note that this new representation is still maximal as a representation in the rank one Lie group ${\rm Aut}(\mathcal Z)$. Indeed, the maximality of the initial representation means that 
$$
\int_X f^\star \o_\Y\wedge c_1(\varphi^\star K_{X_\ca})^{m-1}=\frac{q}{m+1}\, \int_X c_1(\varphi^\star K_{X_\ca})^{m},
$$
but if $\o_{\mathcal Z}$ denotes the complex hyperbolic metric on $\mathcal Z$ normalized to have constant holomorphic sectional curvature $-1$, 
$\o_{\mathcal Z}=\frac{1}{q}\,\o_\Y|_{\mathcal Z}$, so that 
$$
\int_X f^\star \o_{\mathcal Z}\wedge c_1(\varphi^\star K_{X_\ca})^{m-1}=\frac{1}{m+1}\, \int_X c_1(\varphi^\star K_{X_\ca})^{m},
$$
which is exactly the maximality of $\rho:\pi_1(X)\fd{\rm Aut}(\mathcal Z)$. 

If $X_\ca$ were indeed smooth and $\rho(\G)$ discrete and torsion free, the adjunction formula would immediately implies that $K_{X_\ca}=g^\star K_Z+R$ where $Z=\rho(\G)\backslash{\mathcal Z}$ and $R$ is the ramification divisor. In fact, neither the (possible) singularities of $X_\ca$  nor the non discreteness of $\rho(\G)$ are an obstruction for such a formula to hold. Indeed, the line bundle $g^\star K_{\mathcal Z}$ goes down on $X_\ca$ where it defines a Cartier divisor. On the smooth part $X^0_\ca$ of $X_\ca$, thanks to its equivariance, the Jacobian of $g$ defines a (nonzero) section of $K_{X^0_\ca}\otimes g^\star K_{\mathcal Z}^{-1}$ and we denote by $R^0$ its zero divisor. Then, on $X^0_\ca$, we have the relation $K_{X^0_\ca}=(g^\star K_{\mathcal Z})_{|X^0_\ca}+R^0$. Recall that $X_\ca\backslash X^0_\ca$ is at least 2-codimensional, and since $K_{X_\ca}$ and $g^\star K_{\mathcal Z}$ are $\Q$-Cartier, if $R$ denotes the compactification of the divisor $R^0$, we have $K_{X_\ca}=g^\star K_{\mathcal Z}+R$. 

By the maximality of the Toledo invariant, 
$$K_{X_\ca}^m=f^\star K_{\mathcal Z}\cdot(\varphi^\star K_{X_\ca})^{m-1}=g^\star K_{\mathcal Z}\cdot K_{X_\ca}^{m-1}=K_{X_\ca}^m-R\cdot K_{X_\ca}^{m-1}$$
and the ampleness of $K_{X_\ca}$ implies that $R=0$ because $R$ is effective. We conclude that in fact $K_{X_\ca}=g^\star K_{\mathcal Z}$ and $g$ is a local biholomorphism on $\widetilde X_\ca\backslash\widetilde A$.

We then remark that each point $x\in X_\ca$ is isolated in the fiber $g^{-1}(g(x))$. Again, if $\rho(\G)$ is discrete and torsion free we are done because a positive dimensional irreducible component of $g^{-1}(g(x))$ then defines a $s$-dimensional subvariety $S\subset X_\ca$ (for some $s>0$) on which $g^\star K_{\mathcal Z}$ is trivial and $K_{X_\ca}^{s}\cdot S=(g^\star K_{\mathcal Z})^s\cdot S=0$ contradicts the ampleness of $K_{X_\ca}$.
In the general case,  we consider the set $\{x\in\tilde X_\ca\, ,\, {\rm dim}_x g^{-1}(g(x))>0\}$. This set is analytic (see~\cite{Fi}, Theorem 3.6) and $\pi_1(X_\ca)$-invariant thus it defines an analytic subset of $X_\ca$. If it is nonempty and $S$ is a $s$-dimensional irreducible component of it, we also have $K_{X_\ca}^{s}\cdot S=0$.

Now, we can apply Proposition~3.1.2 in \cite{GR}: for any $x\in\widetilde X_\ca$, there exists open neighborhoods $U$ of $x$ in $\widetilde X_\ca$ and $V$ of $g(x)$ in ${\mathcal Z}$ with $g(U)\subset V$, such that the induced map $g_{U,V}:U\fd V$ is a (proper) finite holomorphic map (that is a finite branched covering) and such that $g_{U,V}^{-1}(g_{U,V}(x))=\{x\}$. In fact, it follows from the purity of ramification locus (see the very beginning of \cite{Na}) that for any $x'\in U$, ${\cal O}_{g(x'),V}$ and ${\cal O}_{x',U}$ are isomorphic. Indeed, the ramification locus of $g$, if non empty, should be of pure codimension one and hence would intersect the smooth part of $U$. But we saw above that $g$ is a local biholomorphism on $\widetilde X_\ca\backslash\widetilde A$. Therefore, $U$ is (maybe only locally) biholomorphic to $V$, and $X_\ca$ is smooth.

\medskip

Finally, the map $g:\tilde X_\ca\fd\Bm$ is a local biholomorphism and thus $X_\ca$ can be endowed with a metric of constant holomorphic sectional curvature $-1$, so that it is indeed uniformized by the ball $\Bm$. It follows immediately that $g$ is proper and that $\rho$ is discrete and faithful.
\end{demo}

\end{document}